\documentclass{amsart}
\usepackage{amssymb}
\usepackage[author-year]{amsrefs}
\title[Khintchine's LIL and Extremes]{ An Extreme-Value Analysis of the
		LIL for Brownian Motion}
   \thanks{The research of D.\@ Kh.\@ was
   supported by a grant from the National Science Foundation}
   \author[D.\@ Khoshnevisan]{Davar Khoshnevisan}
   \address{Department\@ of Mathematics, The University of Utah,
      155 S.\@ 1400 E.\@ Salt Lake City, UT 84112--0090}
   \email{davar@math.utah.edu}
   \urladdr{http://www.math.utah.edu/\~{}davar}
   \author[D.\@ A.\@ Levin]{David A. Levin}
   \address{Department\@ of Mathematics, The University of Utah,
      155 S.\@ 1400 E.\@ Salt Lake City, UT 84112--0090}
   \email{levin@math.utah.edu}
   \urladdr{http://www.math.utah.edu/\~{}levin}
\author[Z.\@ Shi]{Zhan Shi}
\address{Laboratoire de Probabilite\'es , 
   Universit\'{e} Paris VI 4 place Jussieu
   F-75252 Paris Cedex 05 France}
\email{zhan@proba.jussieu.fr}
\urladdr{http://www.proba.jussieu.fr/pageperso/zhan}
\date{November 16, 2004}

\theoremstyle{plain}{
\newtheorem{theorem}{Theorem}[section]}
\theoremstyle{plain}{
\newtheorem{proposition}[theorem]{Proposition}}
\theoremstyle{plain}{
   \newtheorem{lemma}[theorem]{Lemma}}
\theoremstyle{plain}{
   }
\theoremstyle{definition}{
   }
\theoremstyle{definition}{
   }
\theoremstyle{remark}{
   \newtheorem{remark}[theorem]{Remark}}
\usepackage{amssymb,euscript}
\usepackage[author-year]{amsrefs}
\numberwithin{equation}{section}

\newcommand{\e}{\varepsilon}

\renewcommand{\P}{\mathrm{P}}
\newcommand{\E}{\mathrm{E}}
\newcommand{\R}{\mathbf{R}}

\renewcommand{\phi}{\varphi}
\newcommand{\itlog}{\mathcal{L}}

\begin{document}
\begin{abstract}
	We present an extreme-value analysis
	of the classical law of the iterated
	logarithm (LIL) for Brownian motion. Our result
	can be viewed as a new improvement to the
	LIL.
\end{abstract}
\maketitle
\section{Introduction}

Let $\{B(t)\}_{t\ge 0}$ be a standard Brownian motion.
The law of the iterated logarithm (LIL)
of~\ocite{khintchine}
states that 
$\limsup_{t\to\infty} (2t\ln\ln t)^{-1/2}B(t)=1$ a.s.
Equivalently,
\begin{equation}\label{lil}
	\text{With probability one,}
	\quad \sup_{s\ge t} \frac{B(s)}{\sqrt{2s\ln\ln s}} \to 1
	\quad\text{as $t\to\infty$} \,.
\end{equation}
The goal of this note is to determine the rate
at which this convergence occurs. 

We consider the extreme-value distribution
function \cite{resnick}*{p.\@ 38},
\begin{equation}
	\Lambda (x) = \exp\left(-e^{-x}
	\right)\qquad {}^\forall x\in\R.
\end{equation}
Also, we place $\mathcal{L}_k x$ or 
$\mathcal{L}_k(x)$ in favor of
the $k$-fold, iterated, natural logarithm,
$(\ln\cdots\ln) (x)$ ($k$ times). Then, our main
result is as follows:

\begin{theorem}\label{main}
   For all $x\in\R$,
   \begin{align}
	\lim_{t\to\infty} \P\left\{
		2 \mathcal{L}_2 t\left(  \sup_{s\ge t}
		\frac{B(s)}{\sqrt{2s \mathcal{L}_2 s}} -1 \right)
		-\frac32\mathcal{L}_3 t + \mathcal{L}_4 t 
		+ \ln\left(\frac{3}{\sqrt{2}}\right)
		\le x\right\} & = \Lambda(x),\label{Eq:Main1}\\
	\lim_{t\to\infty} \P\left\{
		2 \mathcal{L}_2 t\left(  \sup_{s\ge t}
		\frac{|B(s)|}{\sqrt{2s \mathcal{L}_2 s}} -1 \right)
		-\frac32\mathcal{L}_3 t + \mathcal{L}_4 t 
		+ \ln\left( \frac{3}{2\sqrt{2}}\right) \le x\right\} 
		& = \Lambda(x).\label{Eq:Abs}
   \end{align}
\end{theorem}

The preceding is accompanied by
the following strong law:

\begin{theorem}\label{liminf}
	With probability one,
	\begin{equation}
		\lim_{t\to\infty}
		\frac{\mathcal{L}_2 t}{\mathcal{L}_3 t} \left( 
		\sup_{s\ge t} \frac{B(s)}{\sqrt{2s 
		\mathcal{L}_2 s}} -1 \right) 
		= \frac34.
	\end{equation}
\end{theorem}
This should be 
compared with the following consequence
of the theorem of~\ocite{erdos}:
\begin{equation}\label{UC}
	\limsup_{t\to\infty}
	\frac{\mathcal{L}_2 t}{\mathcal{L}_3 t} \left( 
	\sup_{s\ge t}
	\frac{B(s)}{\sqrt{2s
	\mathcal{L}_2 s}} -1 \right) 
	=\frac34\qquad\text{a.s.}
\end{equation}
[Erd\H{o}s's theorem is stated for Bernoulli walks,
but applies equally well---and for the same reasons---to Brownian motion.]

Theorem~\ref{main} is derived by analyzing
the excursions of the Ornstein--Uhlenbeck process,
\begin{equation}
   X(t) = e^{-t/2} B(e^t)\qquad t\ge 0.
\end{equation}
Our method is influenced by the
ideas of~\ocite{motoo}.

\section{Proof of Theorem~\ref{main}}

An application of It\^o's formula shows
us that the process $X$ satisfies the s.d.e,
\begin{equation}\label{SI}
	X(t) = X(0) + \int_1^{\exp(t)} 
	\frac{1}{\sqrt s}\, dB(s) - \frac12 \int_0^t X(s)\, ds.
\end{equation}
The stochastic integral in \eqref{SI}
has quadratic
variation $\int_1^{\exp(t)} s^{-1}\, ds=t$.
Therefore, this  stochastic integral defines
a Brownian motion. Call the said Brownian motion
$W$, to see that
$X$ satisfies the s.d.e.,
\begin{equation}\label{sde}
	dX = dW - \frac12 X\, dt.
\end{equation}
In particular, the quadratic variation of $X$ at
time $t$ is $t$. This means that the semi-martingale
local times of $X$ are occupation densities
\cite{ry}*{Chapter VI}.
In particular, if $\{L_t^0(X)\}_{t\ge 0}$ 
denotes the local time
of $X$ at zero, then
\begin{equation}
	L_t^0(X) = \lim_{\e\to 0} \frac{1}{2\e} \int_0^t
	\mathbf{1}_{\{ |X(s)|\le \e\}}\, ds
	\qquad\text{a.s.}
\end{equation}
See~\ocite{ry}*{Corollary 1.6, p.\@ 224}.

Let $\{\tau(t)\}_{t\ge 0}$ denote the right-continuous
inverse-process to $L^0(X)$. By the ergodic theorem,
$\tau(t)/t$ a.s.\@ converges as $t$ diverges.
In fact,
\begin{equation}\label{tau/t}
	\lim_{t\to\infty}\frac{\tau(t)}{t}=
	\sqrt{2\pi}\qquad\text{a.s.}
\end{equation}
To  this, note that $\tau(t)/t\sim t/L_t^0(X)$
a.s.
But another application of 
the ergodic theorem implies that
$L_t^0(X)\sim \E[ L_t^0(X)]$ a.s. 
The assertion \eqref{tau/t}
then follows from the fact that
$\E [L_t^0(X)] =t/\sqrt{2\pi}$~\cite{hk}*{Lemma 3.2}.

Define
\begin{equation}
	\e_t = \sup_{s\ge \tau(t)}
	\frac{X(s)}{\sqrt{2\ln s}}\qquad
	{}^\forall t\ge e.
\end{equation}

\begin{lemma}\label{theta}
	Almost surely,
	\begin{equation}
		\left| \e_n - \sup_{j\ge n}
		\frac{M_j}{\sqrt{2\ln j}}
		\right| = O\left(\frac{1}{\ln n}\cdot
		\sqrt{\frac{\mathcal{L}_2 n}{n}} \right)
		\qquad(n\to\infty),
	\end{equation}
	where	
	\begin{equation}\label{Mj}
		M_j = \sup_{s\in [\tau(j),\tau(j+1)]}X(s)
		\qquad {}^\forall j\ge 1.
	\end{equation}
\end{lemma}

\begin{proof}
	According to \eqref{tau/t},
	\begin{equation}
		\sup_{s\in[\tau(j),\tau(j+1)]}
		\left| \frac{1}{\sqrt{\ln\tau(j)}}
		-\frac{1}{\sqrt{\ln s}}\right|
		\sim \frac{1}{2\ln j}\cdot
		\sqrt{\frac{\mathcal{L}_2 j}{j\ln j}}
		\qquad(j\to\infty).
	\end{equation}
	On the other hand, according
	to~\eqref{lil} and~\eqref{tau/t}, almost surely,
	\begin{equation}
		M_j = O\left( \sqrt{\ln \tau(j+1)}\right)
		=O\left(\sqrt{\ln j}\right)\qquad
		(j\to\infty).
	\end{equation}
	The lemma follows from a little algebra.
\end{proof}

Lemma~\ref{theta}, and monotonicity,
together prove that Theorem~\ref{main} is equivalent to
the following: For all $x\in\R$,
\begin{equation}\label{goal}
	\lim_{n\to\infty}  \P\left\{
	2\ln n\left( \sup_{j\ge n} \frac{M_j}{\sqrt{2\ln j}}
	-1 \right) -\frac32 \mathcal{L}_2 n +
	\mathcal{L}_3 n + \ln\left(
	\frac{3}{\sqrt{2}} \right) \le x\right\}\\
	 = \Lambda (x).
\end{equation}
We can derive this because: (i)
By the
strong Markov property of the OU process
$X$, $M_1,M_2,\ldots$ is an i.i.d.\@ sequence;
and (ii) the distribution of $M_1$ can be found
by a combining a little bit of
stochastic calculus with an iota of excursion theory.
In fact, one has a slightly more general result
for It\^o diffusions (i.e., diffusions
that solve smooth s.d.e.'s)
at no extra cost.

\begin{proposition}
	Let $\{Z_t\}_{t\ge 0}$ denote the regular It\^o diffusion
	on $(-\infty,\infty)$ which solves the s.d.e.
	\begin{equation} \label{Eq:SDE}
		dZ_t = \sigma(Z_t) \, dw_t + a(Z_t)dt \,,
	\end{equation}
	where $\sigma,a \in \mathcal{C}^\infty(\R)$,
	$\sigma$ is bounded away from zero, and  $\{w_t\}_{t\ge 0}$
	is a Brownian motion.
	Write $\{\theta_t\}_{t\ge 0}$ for
	the inverse local-time of $\{Z_t\}_{t \geq 0}$ at zero.
	Then for all $\lambda>0$,
	\begin{equation}
		\P\left\{ \left. \sup_{t\in[0,\theta_1]}
		Z_t \le\lambda \ \right|\ Z_0=0 \right\} 
		 = \exp\left(
		-\frac{ f'(0)}{2 \left\{ f(\lambda)-f(0)\right\}}
		\right).
	\end{equation}

\end{proposition}

\begin{proof}
	The scale function of a diffusion is defined
	only up to an affine transformation. Therefore,
	we can assume, without loss of generality,
	that $f'(0)=1$ and $f(0)=0$; else, we choose the
	scale function $x\mapsto \{ f(x)-f(0) \}/f'(0)$
	instead.  Explicitly, $\{Z_t\}_{t\ge 0}$ has the
	scale function \cite{ry}*{Exercise VII.3.20}.
	\begin{equation} \label{Eq:ScaleFunction}
	  f(x) = \int_0^x \exp\left( - 2 \int_0^y
	  \frac{a(u)}{\sigma^2(u)} du \right)\, dy.
	\end{equation}

	Owing to It\^o's formula, $N_t:=f(Z_t)$ satisfies
	\begin{equation} \label{Eq:DN}
		dN_t = f'(Z_t)\, \sigma(Z_t)\, dw_t = 
		f'\left( f^{-1}(N_t) \right)\,
		\sigma\left( f^{-1}(N_t) \right)\, dw_t \,,
	\end{equation}
	and so is a local martingale.
	According to the
	Dambis, Dubins--Schwarz representation theorem
	\cite{ry}*{Theorem V.1.6, p.\@ 181}, 
	there exists a Brownian motion
	$\{b(t)\}_{t\ge 0}$ such that
	\begin{equation} \label{Na}
	\begin{split}
		N_t & = b(\alpha_t) \,, \quad\text{where} \\
		\alpha_t & = \alpha(t) = \langle N \rangle_t = 
		\int_0^t \left[ f'\left(f^{-1}
		(N_r) \right) \right]^2\,
		\sigma^2\left(f^{-1}(N_r)\right)\, dr\qquad{}^\forall t\ge 0.
	\end{split}
	\end{equation}
	The process $N$ is manifestly a diffusion;
	therefore, it has continuous local-time processes
	$\{L_t^x(N)\}_{t\ge 0, x\in\R}$ which satisfy the
	occupation density formula~\cite{ry}*{Corollary VI.1.6, p.\@ 224
	and time-change}.
%	In particular,
%	\begin{equation}
%		\alpha_t = \int_{-\infty}^\infty
%		\left[ f'\left( f^{-1}(x) \right)\right]^2
%		\, \sigma^2\left( f^{-1}(x) \right)
%		L_t^x(N)\, dx\qquad{}^\forall t\ge 0.
%	\end{equation}
	By \eqref{Eq:ScaleFunction}, $f' > 0$,
%	Because $f$ is strictly increasing and $\mathcal{C}^2$,
	and because $\sigma$ is bounded away from zero,
	$\sigma^2 f'> 0$. Therefore, the inverse process
	$\{\alpha^{-1}(t)\}_{t\ge 0}$ exists a.s., and is uniquely
	defined by $\alpha(\alpha^{-1}(t))=t$ for all $t\ge 0$.
	
	Let $\{L_t^x(b)\}_{t\ge 0, x\in\R}$ denote the
	local-time processes of the Brownian motion
	$b$. It is well known \cite{RW:DMM2}*{Theorem V.49.1} that a.s.
	\begin{equation} \label{Eq:LTEquality}
		\{L_t^x(b)\}_{t \geq 0, x \in \R}
		= \{L_{\alpha(t)}^x(N)\}_{t \geq 0, x \in \R} \,;
	\end{equation}
	For completeness, we include a brief argument here.
	Thanks to the occupation density formula,
	\begin{equation}\begin{split}
		\int_{-\infty}^\infty h(x) L_t^x(b)\, dx &=
		\int_0^t h(b(s))\, ds =
		\int_0^t h(N_{\alpha^{-1}(s)})\, ds,
	\end{split}\end{equation}
	valid for all Borel-measurable functions $h:\R\to\R$.
	See~\eqref{Na} for the last equality. We can change variables
	[$v=\alpha^{-1}(s)$], and use the definition of 
	$\alpha_t$ in \eqref{Na}, to note that
	\begin{equation}\begin{split}
		\int_{-\infty}^\infty h(x) L_t^x(b)\, dx &=
			\int_0^{\alpha^{-1}(t)} h(N_v)\, d\alpha_v\\
%		&= \int_0^{\alpha^{-1}(t)} h(N_v)
%			\left[ f'\left(f^{-1}(N_v)\right)\right]^2
%			\, \sigma^2\left( f^{-1}(N_v) \right) \, dv\\
		&= \int_0^{\alpha^{-1}(t)} h(N_v) d\langle N \rangle_v
			= \int_{-\infty}^\infty h(x)
			L_{\alpha^{-1}(t)}^x(N)\, dx.
	\end{split}\end{equation}
	This establishes \eqref{Eq:LTEquality}.	
%	Therefore, with probability one,
%	\begin{equation}
%		L_t^x (b) = 
%		L_{\alpha^{-1}(t)}^x(N)\qquad
%		{}^\forall t\ge 0,~ x\in\R.
%	\end{equation}
	In particular, with probability one,
	\begin{equation}\label{LNLb}
		 L_t^0(N) = L_{\alpha(t)}^0(b)\qquad
		{}^\forall t\ge 0.
	\end{equation}
	By \eqref{Eq:SDE} and \eqref{Eq:DN},
	$d\langle Z \rangle_t = 
	\left[ f'\left(f^{-1}(N_t)\right) \right]^{-2} d \langle N \rangle_t$,
	so if $L(Z)$ denotes the local times of
	$Z$, then almost surely,
	\begin{equation} \begin{split}
		\int_{-\infty}^\infty h(x) L_t^x(Z)\, dx
			& = \int_0^t h(Z_r)\, d\langle Z \rangle_r \\
			& = \int_0^t \frac{h\left( f^{-1}(N_r) 
				\right)}{\left[ f'\left(f^{-1}(N_r)\right) 
			 	\right]^2} \, d\langle N \rangle_r \\
			& = \int_{-\infty}^\infty h\left( f^{-1}(y) \right)
			\frac{ L_t^y(N) }{\left[ f'\left(f^{-1}(y)\right) \right]}
			\, d\left[ f^{-1}(y) \right] \\
			& = \int_{-\infty}^\infty h(x) \frac{ L_t^{f(x)}(N) }{f'(x)}
			\, dx\quad [x=f^{-1}(y)].
	\end{split} \end{equation}
	This proves that a.s.,
	$f'(x) \cdot L_t^x(Z) = L_t^{f(x)}(N)$. In particular, $L_t^0(Z) = L_t^0(N)$
	for all $t\ge 0$, a.s. Thanks to \eqref{LNLb},
	we have proved the following: Almost surely,
	\begin{equation}\label{LNLz}
		L_t^0(Z) = L_{\alpha(t)}^0(b)\qquad{}^\forall t\ge 0.
	\end{equation}
	
	Define
	\begin{equation}
		\phi_t = \inf\left\{ s>0:\
		L_s^0(b) > t\right\}
		\qquad{}^\forall t\ge 0.
	\end{equation}
	Then thanks to \eqref{LNLz},
	\begin{equation}\label{tau}
		\phi_t = \alpha\left(\theta_{t}\right)\qquad {}^\forall t\ge 0.
	\end{equation}
	Thus,
	\begin{equation} \label{Eq:SupZ} \begin{split}
		\P\left\{ \left. \sup_{s\in[
			0,\theta_1]} Z_s \le \lambda \ \right|\
			Z_0=0 \right\} &= \P\left\{ \left. \sup_{s\in[
			0,\theta_1]} N_s \le f(\lambda) \ \right|\
			N_0=0 \right\}\\
		&= \P\left\{ \left. \sup_{s\in\left[
			0,\phi_1 \right]} b_s \le f(\lambda) \ \right|\
			b_0=0 \right\}.
	\end{split}\end{equation}
	The last identity follows from \eqref{Na}, and the fact
	that $\alpha$ and $\alpha^{-1}$ are both continuous
	and strictly increasing a.s. 
	
	Define $\mathcal{N}_\beta$
	to be the total number of excursion of the 
	Brownian motion $b$ that exceed
	$\beta$ by local-time $1$. Then,
	\begin{equation}\begin{split} \label{Eq:SupBM}
		\P\left\{ \left. 
		\sup_{s\in\left[
			0,\phi_1 \right]} b_s
		\le f(\lambda) \ \right|\ 
			b_0=0 \right\} &= 
			\P\left\{ \left.
			\mathcal{N}_{f(\lambda)} =0\
			\right|\ b_0=0 \right\}\\
		&= \exp\left\{- \E
			\left[ \left. \mathcal{N}_{f(\lambda)}
			\, \right|\, b_0=0 \right] \right\},
	\end{split}\end{equation}
	because $\mathcal{N}_\beta$
	is a Poisson random variable~\cite{ito}.
	According to
	Proposition 3.6 of~\ocite{ry}*{p.\@ 492},
	$\E[\mathcal{N}_\beta \,|\, b_0=0]=(2\beta)^{-1}$ for all
	$\beta>0$. [See also \ocite{ry}*{Exercise XII.4.11}.]
	The result follows.
\end{proof}

\begin{remark} \label{Rmk:Abs1}
	Also, the following equality holds:
	\begin{equation} \label{Eq:Abs2}
			\P\left\{ \left. \sup_{t\in[0,\theta_1]}
			|Z_t| \le\lambda \ \right|\ Z_0=0 \right\} 
			 = \exp\left(
			-\frac{ f'(0)}{ f(\lambda)-f(0) }
			\right).
	\end{equation}
	This follows as above after noting that
	$f(-x) = -f(x)$, and that
	$\E[{\mathcal N}'_\beta \, | \, b_0 = 0] = \beta^{-1}$,
	where
	${\mathcal N}'_\beta$ denotes
	the number of excursions of the Brownian motion $b$ that exceed
	$\beta$ in absolute value by local-time $1$.
\end{remark}

\begin{proof}[Proof of Theorem \ref{main}]
	If we apply the preceding computation to the diffusion
	$X$ itself, then we find that
	$\P\{M_1\le\lambda\}=\exp\{ -1/[2S(\lambda)]\}$,
	where $S$ is the scale function of $X$ which satisfies
	$S'(0)=1$ and $S(0)=0$. According
	to \eqref{sde} 
	and \eqref{Eq:ScaleFunction}, $S(x)=\int_0^x \exp(y^2/2)\, dy$.
	Note that $S(x)\sim x^{-1}\exp(x^2/2)$ as $x\to\infty$. 
	
	Let $\{\beta_n(x)\}_{n=1}^\infty$ be 
	a sequence which, for $x$ fixed,
	satisfies $\beta_n(x) \rightarrow \infty$ as $n \rightarrow \infty$.
	We assume, in addition,
	that $\alpha_n(x) := \beta_n(x)/\ln n$ goes to zero
	as $n \rightarrow \infty$.
	We will suppress $x$ in the notation and write $\alpha_n$ and
	$\beta_n$ for $\alpha_n(x)$ and $\beta_n(x)$, respectively.
	
	A little calculus shows that if $\alpha_n(x) > 0$, then
	\begin{equation}\label{prod}\begin{split}
	   \P\left\{ \sup_{j\ge n} \frac{M_j}{\sqrt{
		2\ln j}} \le 1+\frac{\alpha_n}{2} \right\}
		& = \prod_{j=n}^\infty \exp\left( -
		\frac{1}{2 S\left( (1 + \alpha_n/2)\sqrt{2 \ln j}
		\right)} \right)\\
	   & = \exp\left( -\frac{(1 + \alpha_n/2)}{\sqrt 2}\sum_{j=n}^\infty
		\frac{ \sqrt{\ln j}}{j^{(1+\alpha_n/2)^2}} \right)\\
	   & = \exp\left( -\frac{[1+o(1)] 
	   \left(1 + \alpha_n/2 \right) }{ \sqrt 2}\int_n^\infty
	        \frac{\sqrt{\ln x} }{x^{(1+ \alpha_n/2)^2} } dx\right)\\
	   & = \exp\left( - \frac{ [1 + o(1)] 
	      \left( 1 + \alpha_n/2 \right) \sqrt{\ln n} }{
	      \sqrt{2} \alpha_n \left( 1 + \alpha_n/4 \right)
	   	\alpha_n n^{\alpha_n\left( 1 + \alpha_n/4 \right)} }
		\right) \\
	   & = \exp\left( - \left[ 1 + o(1)\right] 
	   	q_n(x) \frac{ (\ln n)^{3/2} 
		e^{-\beta_n} }{ \sqrt{2} \beta_n} \right).
	\end{split}\end{equation}
	Here, 
	\begin{equation}
		q_n(x) = \left(\frac{2\ln n + \beta_n(x)}{
		2\ln n + \beta_n(x)/2} \right) \exp
		\left(- \frac{\beta^2_n(x)}{4\ln n} \right) \,.
		\end{equation}
	If $\alpha_n(x) \leq 0$, then the probability on the right-hand side of
	\eqref{prod} is $0$.
	Define
	\begin{equation} \label{Eq:Beta}
		\phi_n(x) := \frac32 \mathcal{L}_2 n - 
		\mathcal{L}_3 n -\ln\left(3/\sqrt{2}\right) +x,
	\end{equation}
	and set $\beta_n(x)$ in \eqref{prod} equal 
	to $\phi_n(x)$. This yields
	\begin{equation} \label{Eq:MBeta}
		\P\left\{ \sup_{j\ge n} \frac{M_j}{\sqrt{
		2\ln j}} \le 1+\frac{\phi_n(x)}{2\ln n} \right\} 
	 	= 
		\begin{cases} 
			\exp\left( -[1 + o(1)] c_n(x) 
				e^{-x} \right) & \text{if } \phi_n(x) > 0,\\
			0 & \text{if } \phi_n(x) \leq 0,
		\end{cases}
	\end{equation}
	where
	\begin{equation} \label{Eq:c}
		c_n(x) =
		\left( \frac{2\ln n + \phi_n(x)}
		{ 2\ln n + \phi_n(x)/2 } \right) 
		\exp\left( - \frac{\phi_n^2}{4 \ln n} \right) \left[
		1 + \frac{-\itlog_3n - \ln(3/\sqrt{2}) + x}{\frac{3}{2} 
		\itlog_2 n} \right]^{-1}.
	\end{equation}
	Note that the little-$o$ in \eqref{Eq:MBeta} is uniform in $x$.
	If $x \in \R$ is fixed, letting $n \rightarrow \infty$ in
	\eqref{Eq:MBeta} shows that
	\begin{equation}
	\lim_{n \rightarrow \infty} 
		\P\left\{ 2 \ln n\left( \sup_{j \geq n} \frac{M_j}{\sqrt{2 \ln j}}
		- 1 \right) \leq \phi_n(x) \right\} =
		\exp\left(-e^{-x}\right) \,.
	\end{equation}
	This proves~\eqref{goal},
	whence equation \eqref{Eq:Main1} of Theorem~\ref{main} follows.
	
	Using \eqref{Eq:Abs2}, we obtain also
	\begin{equation} \label{Eq:AbsMBeta}
	   \P\left\{ \sup_{j\ge n} \frac{|M_j|}{\sqrt{
		2\ln j}} \le 1+ \frac{\alpha_n}{2} \right\}
	   = \exp\left( -[1+o(1)]\sqrt{2}
		\cdot \frac{e^{-\beta_n}\left(\ln n\right)^{3/2}}{\beta_n}
		\right) \,.
	\end{equation}
	Let $\beta_n(x) = \frac32 \mathcal{L}_2 n - 
	\mathcal{L}_3 n -\ln\left(\frac{3}{2\sqrt{2}}\right) +x$ in
	\eqref{Eq:AbsMBeta} to establish~\eqref{Eq:Abs}.
\end{proof}

\section{Proof of Theorem~\ref{liminf}}
In light of~\eqref{UC} it suffices to prove
that
\begin{equation}
	\liminf_{t\to\infty}
	\frac{\mathcal{L}_2 t}{\mathcal{L}_3 t} \left( 
	\sup_{s\ge t} \frac{B(s)}{\sqrt{2s 
	\mathcal{L}_2 s}} -1 \right) 
	\ge \frac34\qquad\text{a.s.}
\end{equation}

We aim to prove that almost surely,
\begin{equation}\label{lb}
	\sup_{j\ge n} \frac{M_j}{\sqrt{2\ln j}}
	> \sqrt{ 1+ 
	c\frac{\mathcal{L}_2 n}{\ln n}}
	\quad\text{eventually a.s.\@ if $c<\frac32$}
\end{equation}
Theorem~\ref{liminf} follows from this by the similar
reasons that yielded Theorem~\ref{main} from~\eqref{goal}.
But \eqref{lb} follows from \eqref{prod}:
\begin{equation}
	\P\left\{ \sup_{j\ge n} \frac{M_j}{\sqrt{2\ln j}}
	\le \sqrt{ 1+ 
	c\frac{\mathcal{L}_2 n}{
	\ln n}} \right\}
	= \exp\left\{ - \frac{1+o(1)}{2c\sqrt 2}\cdot
	\frac{(\ln n)^{(3/2)-c}}{\mathcal{L}_2 n}\right\}.
\end{equation}
Replace $n$ by $\rho^n$ where $\rho>1$ is fixed.
We find that if $c<(3/2)$ then the probabilities
sum in $n$. Thus, by the Borel--Cantelli lemma,
for all $\rho>1$ and $c<(3/2)$ fixed,
\begin{equation}
	\sup_{j\ge \rho^n} \frac{M_j}{\sqrt{2\ln j}}
	> \sqrt{ 1+ 
	c\frac{\mathcal{L}_2 (\rho^n)}{\ln (\rho^n)}}
	\quad\text{eventually a.s.}
\end{equation}
Equation~\eqref{lb} follows from this and
monotonicity.

\section{An Expectation Bound}

We can use our results to improve on the bounds
of Dobric and Marano \ycite{dobricmarano} for
the rate of convergence of 
$\E[\sup_{s \geq t}
 B_s (2 s\itlog_2 s )^{-1/2} ]$ 
to $1$.   
\begin{proposition}
	As $t\to\infty$,
	\begin{equation}
		\E\left[\sup_{s\ge t} \frac{B(s)}{
		\sqrt{2s\itlog_2 s}}\right]
		= 1+\frac34 \frac{\itlog_3 t}{\itlog_2 t}
		-\frac12 \frac{\itlog_4 t}{\itlog_2 t} 
		+\frac12
		\frac{\gamma- \ln\left(3/\sqrt{2}\right)}{\itlog_2 t}
		+o\left(\frac{1}{\itlog_2 t}\right),
	\end{equation}
	where $\gamma\approx 0.5772$ denotes Euler's constant.
\end{proposition}

\begin{proof}
	Define
	\begin{equation}
		U_n := 2 \ln n \left( \sup_{j \geq n} \frac{M_j}{\sqrt{2\ln j}} - 1
		\right) - \frac{3}{2} \itlog_2 n + \itlog_3 n + \ln\left(3/\sqrt{2}\right)\,. 
	\end{equation}
	We have shown that $U_n$ converges weakly to $\Lambda$.
	We now establish that $\sup_n \E\left[ U_n^2
	\right] < \infty$. This implies uniform integrability,
	whence we can deduce that $\E[U_n] 
	\rightarrow \int x\, d\Lambda(x)$.
	
	Let $\phi_n(x)$ be as defined in \eqref{Eq:Beta}, and
	note that
	$x_n^\star := (3/2) \itlog_2n -
	\itlog_3 n - \ln(3/\sqrt{2})$ solves $\phi_n(-x_n^\star) = 0$. 
	  
	Recalling the definition of $c_n(x)$ in \eqref{Eq:c} and
	rewriting \eqref{Eq:MBeta} shows that for $x > -x_n^\star$
	\begin{equation}
		P\left\{ U_n \leq x \right\} = 
		\exp\left( -[1+o(1)]c_n(x) e^{-x} \right).
	\end{equation}
	Consequently, for $n$ large enough,
	\begin{equation}
		\int_0^\infty x\P\left\{ U_n \leq - x\right\}  \leq
		\int_0^{x_n^\star} x \exp\left( - 
		\frac{1}{2} c_n(-x) e^x\right)\, dx.
	\end{equation}
	For $n$ sufficiently large and $0 \leq x < x_n^\star$,
	$c_n(-x) \geq e^{-2}/(3/2) \geq (1/12)$. 
	Thus, for $n$ sufficiently large,
	\begin{equation}
	\int_0^\infty x\P\left\{ U_n \leq - x\right\}  \leq \int_0^{x_n^\star}
	x e^{- \frac{1}{24} e^{x} } dx \leq \int_0^\infty x e^{- \frac{1}{24}
	e^{x} } dx < \infty \,.
	\end{equation}
	Also for $n$ sufficiently large,
	\begin{equation}
		\int_0^\infty x \P\left\{ U_n > x \right\} \leq
		\int_0^\infty x \left( 1 - e^{- \frac{3}{2} 
		c_n(x) e^{-x}}\right) dx \leq
		\int_0^\infty \frac{3}{2}x c_n(x) e^{-x}\, dx.
	\end{equation}
	We can get the easy bound for $x > 0$, 
	$c_n(x) \leq (3/2 + x)/(1/4) = (6 + 4x)$, yielding
	\begin{equation}
		\int_0^\infty x \P\{ U_n > x\}\, dx 
		\leq \int_0^\infty \frac{3}{2} x\left(6 + 4x \right)
		e^{-x}\, dx < \infty.
	\end{equation}
	Now we can write
	\begin{equation}
		\E\left[ U_n^2 \right] = 2\int_0^\infty x
		\P\left\{ | U_n| > x \right\}\, dx
		= 2\int_0^\infty x\P\left\{ U_n > x\right\}\, dx 
		+ 2\int_0^\infty x\P\left\{
		 U_n < -x \right\}\, dx.
	 \end{equation}
	We have just show that the two terms on the right-hand side are
	bounded uniformly in $n$, which establishes uniform integrability.

	From Lemma \ref{theta}, it follows that
	\begin{equation}
		2\itlog_2 t \left( \E\left[\sup_{s \geq t} 
		\frac{B_s}{\sqrt{2 s\itlog_2 s}}\right] - 1 \right) -
		\frac{3}{2} \itlog_3t + \itlog_4 t + \ln\left(
		\frac{3}{\sqrt{2}}\right) 
		\rightarrow \int_{-\infty}^\infty x\, d\Lambda(x).
	\end{equation}
	It suffices to prove that $\int x\, d\Lambda(x)=\gamma$.
	But this follows because
	\begin{equation}
		\int_{-\infty}^\infty x\, d\left[ e^{-e^{-x}}\right]
		= - \int_{0}^{\infty} \ln(t) e^{-t}\, dt
		= - \frac{d}{dz} \int_{0}^\infty t^{z-1} e^{-t}\, dt
		\Bigg|_{z = 1} 
		= - \frac{\Gamma'(z)}{\Gamma(z)} \Bigg|_{z = 1}.
	\end{equation}
	Cf.\@~\ocite{gamma}*{Equation 6.3.2}. 
\end{proof}

\section{An Application to Random Walks}
Let $X_1,X_2,\ldots$ be i.i.d.\@ random variables
with 
\begin{equation}\label{moments}
	\E[X_1]=0,\ \mathrm{Var}(X_1)=1,\
	\text{ and }\
	\E\left[ X_1^2\itlog_2\left(X_1 \vee e^e\right)\right]<\infty.
\end{equation}
Let $S_n=X_1+\cdots+X_n$ ($n\ge 1$) denote the corresponding
random walk. Then,
according to Theorem 2 
of~\ocite{einmahl}, there exists a probability space
on which one can construct $\{S_n\}_{n=1}^\infty$
together with a Brownian motion $B$ such that
$| S_n - B(n) |^2 = o( 
{n/\itlog_2 n})$ a.s.
On the other hand, by the Borel--Cantelli
lemma, $|B(t)-B(n)|^2 =o({n/\itlog_2 n})$ 
uniformly for all $t\in[n,n+1]$ a.s.
These remarks,
and a few more lines of elementary computations,
together yield the following.

\begin{proposition}
   If~\eqref{moments} holds, then for all $x\in\R$,
   \begin{align}
	\lim_{n\to\infty} \P\left\{
		2 \mathcal{L}_2 n\left(  \sup_{k\ge n}
		\frac{S_k}{\sqrt{2k \mathcal{L}_2 k}} -1 \right)
		-\frac32\mathcal{L}_3 n + \mathcal{L}_4 n 
		+ \ln\left(\frac{3}{\sqrt{2}}\right)
		\le x\right\} & = \Lambda(x),\\%\label{Eq:Main1}\\
	\lim_{n\to\infty} \P\left\{
		2 \mathcal{L}_2 n\left(  \sup_{k\ge n}
		\frac{|S_k|}{\sqrt{2k \mathcal{L}_2 k}} -1 \right)
		-\frac32\mathcal{L}_3 n+ \mathcal{L}_4 n 
		+ \ln\left( \frac{3}{2\sqrt{2}}\right) \le x\right\} 
		& = \Lambda(x).
   \end{align}
\end{proposition}

\begin{remark}
	It would be interesting to know if the preceding 
	remains valid if only
	$\E[X_1]=0$ and $\mathrm{Var}(X_1)=1$. We believe
	the answer to be, ``No.''
\end{remark}

\end{document}